\documentclass[12pt,leqno]{article}
\usepackage{latexsym,amsfonts,epsfig}
\usepackage{amsmath,amssymb,amsthm}
\usepackage[notref, notcite]{}
\usepackage[colorlinks,linkcolor=red,citecolor=blue,urlcolor=blue]{hyperref}
\usepackage[margin=0.5in]{geometry}
\usepackage{graphicx}
\input epsf
\newtheorem{theorem}{Theorem}[section]

\newtheorem{corollary}{Corollary}[section]
\newtheorem{remark}{Remark}[section]

\RequirePackage[colorlinks,citecolor=blue,urlcolor=blue]{hyperref}

\begin{document}

\title{Sharp unifying generalizations of Opial's inequality}
\author{ Chris A.J. Klaassen \\
Korteweg-de Vries Insitute for Mathematics \\
University of Amsterdam}

\maketitle

\date{}



\begin{abstract}
Opial's inequality and its ramifications play an important role in the theory of differential and difference equations.
A sharp unifying generalization of Opial's inequality is presented that contains both its continuous and discrete version.
This generalization based on distribution functions is extended to the case of derivatives of arbitrary order.
This extension optimizes and improves the constant as given in the literature.
The special case of derivatives of second order is studied in more detail.
Two closely related Opial inequalities with a weight function are presented as well.
The associated Wirtinger inequality is studied briefly.
\end{abstract}

\section{Introduction}

More than sixty years ago \cite{Opial2} published an inequality that happened to have an enormous influence in the field of differential and difference equations.
For an early application see \cite{Willett}.
Its complicated proof was simplified in the early 1960s by several authors; see Chapter 1 of \\
\cite{Agarwal}, who present a comprehensive account of the state of the art in this field up to 1995.
A more recent detailed description of the history of Opial's inequality is given in \cite{Saker}.
The situation in 2015 is summarized in \cite{Andric}.
The inequality has been extended and generalized in several directions.
For an exposition of Opial's inequality for time scales we refer to \cite{Osman}.

In the next section we present a very general, unifying version of the original Opial inequality in terms of integration via distribution functions.
Observe that a distribution function may have its support in different time scales, like ${\mathbb R}, {\mathbb N}$ and $q^{\mathbb N}$.
In particular, our probability theoretic inequality incorporates also a discrete version, which is similar to the one presented by \cite{Lasota}.
We also present the Opial inequality for $n$-th derivatives, with the optimal constant.
Our inequality improves the constants as given in the literature, e.g. \cite{Das}; see also \cite{Wang}.
This is done in Section 3.
However, this inequality is not sharp for distributions with atoms, in particular for the discrete case.
An inequality for derivatives of second order that is sharp for these distributions as well, is proved in Section 4.
Furthermore, we consider Opial inequalities with weights in Section 5.
We study the related Wirtinger inequality briefly in Section 6.
Our conclusions are presented in Section 7.

Formula (8) in \cite{Opial1} states that for differentiable functions $x$ with $x(0) = x(h) =0$ the inequality
\begin{equation}\label{A1}
\int_0^h |x(t) x'(t)| \, dt \leq \frac h \pi \int_0^h |x'(t)|^2\, dt, \quad x(0)=x(h)=0,
\end{equation}
holds.
As mentioned in (1)--(4) of \cite{Opial2} Item 257 on page 185 of\\ \cite{Hardy} implies, by the substitution $x(t) = y(\pi t/h)$,
the Wirtinger inequality
\begin{equation}\label{A2}
\int_0^h |x(t)|^2 \, dt \leq \frac{h^2}{\pi^2} \int_0^h |x'(t)|^2\, dt, \quad x(0)=x(h)=0,
\end{equation}
which, by Cauchy-Schwarz, yields (\ref{A1}); cf. (7) and (8) of \cite{Opial1}.
In (5) of \cite{Opial2} inequality (\ref{A1}) has been sharpened to
\begin{equation}\label{A3}
\int_0^h |x(t) x'(t)| \, dt \leq \frac h4 \int_0^h |x'(t)|^2\, dt, \quad x(0)=x(h)=0,
\end{equation}
where $h/4$ is the optimal constant.
His complicated proof has been simplified successively by \cite{Olech}, \cite{Beesack}, \cite{Levinson}, \cite{Mallows} and \cite{Pederson}.

Actually, (\ref{A3}) is a consequence of the more fundamental inequality for differentiable functions $x$ with just $x(0) =0$ (cf. (3) of \cite{Beesack})
\begin{equation}\label{A4}
\int_0^h |x(t) x'(t)| \, dt \leq \frac h2 \int_0^h |x'(t)|^2\, dt,
\end{equation}
where $h/2$ is the optimal constant; see also Remark \ref{R1}.
Often this inequality is also called Opial's inequality.
It is this inequality we will focus on, but consequences such as (\ref{A3}) will be discussed as well.

\section{A probability theoretic generalization of Opial's inequality}

In this section we first generalize Opial's fundamental inequality (\ref{A4}) to integration via distribution functions.
Subsequently, we will discuss a simple extension to an analog of (\ref{A3}).
Opial-type inequalities are, for obvious reasons, presented in terms of a function $x$ and its derivatives.
To formulate our results we start with the derivative, which we call $\psi$, and integrate this function with respect to a distribution function $F$. Here we define $F(x) = P_F(X \leq x)$ at $x \in {\mathbb R}$ as the probability that a random variable under $F$ equals at most $x$.
Note that $F$ is continuous from the right with left hand limits.
To compensate for this asymmetry in the standard definition of a distribution function, we integrate a measurable function $\psi$ with respect to such a distribution function in a special way, namely
\begin{equation}\label{A}
\int_{-\infty}^\infty \psi(y) \left\{ {\bf 1}_{[y<x]} + \tfrac 12 {\bf 1}_{[y=x]} \right\}\, dF(y)
= E_F \left( \psi(Y) \left\{ {\bf 1}_{[Y<X]} + \tfrac 12 {\bf 1}_{[Y=X]} \right\} \mid X=x \right).
\end{equation}
We will also consider
\begin{equation}\label{B}
\int_{-\infty}^\infty \psi(y) \left\{ {\bf 1}_{[y>x]} + \tfrac 12 {\bf 1}_{[y=x]} \right\}\, dF(y)
= E_F \left( \psi(Y) \left\{ {\bf 1}_{[Y>X]} + \tfrac 12 {\bf 1}_{[Y=X]} \right\} \mid X=x \right).
\end{equation}

\begin{theorem}\label{Opial}
{\bf Opial's inequality} \\
Let $X$ and $Y$ be independent and identically distributed (i.i.d.) random variables with distribution function $F$ on $\mathbb R$, and let $\psi : {\mathbb R} \to {\mathbb R}$ be a measurable function. The inequalities
\begin{eqnarray}\label{O12.1.1}
\lefteqn{ E_F \left( \left| E_F \left( \psi(Y) \left\{ {\bf 1}_{[Y<X]} + \tfrac 12 {\bf 1}_{[Y=X]} \right\} \mid X \right) \psi(X) \right| \right) } \\
&& \hspace{5em} \leq E_F \left( |\psi(X)\psi(Y)| \left\{ {\bf 1}_{[Y<X]} + \tfrac 12 {\bf 1}_{[Y=X]} \right\} \right)
\leq \tfrac 12 E_F \left( \psi^2(X)\right) \nonumber
\end{eqnarray}
and
\begin{eqnarray}\label{O12.1.2}
\lefteqn{ E_F \left( \left| E_F \left( \psi(Y) \left\{ {\bf 1}_{[Y>X]} + \tfrac 12 {\bf 1}_{[Y=X]} \right\} \mid X \right) \psi(X) \right| \right) } \\
&& \hspace{5em} \leq E_F \left( |\psi(X)\psi(Y)| \left\{ {\bf 1}_{[Y>X]} + \tfrac 12 {\bf 1}_{[Y=X]} \right\} \right)
\leq \tfrac 12 E_F \left( \psi^2(X)\right) \nonumber
\end{eqnarray}
hold with equalities if and only if $\psi$ is constant $F$-almost everywhere.
\end{theorem}
\noindent
{\bf Proof} \\
Without loss of generality we assume $E_F\left( \psi^2(X) \right) < \infty$.
Since the absolute value of a (conditional) expectation of a random variable is bounded from above by the expectation of the absolute value of the random variable itself, the left-hand side of (\ref{O12.1.1}) is bounded by the middle term
\begin{equation}\label{O12.2}
E_F \left( |\psi(X)\psi(Y)| \left\{ {\bf 1}_{[Y<X]} + \tfrac 12 {\bf 1}_{[Y=X]} \right\} \right)
= E_F \left( |\psi(X)\psi(Y)| \left\{ {\bf 1}_{[Y>X]} + \tfrac 12 {\bf 1}_{[Y=X]} \right\} \right),
\end{equation}
where the equality holds because $X$ and $Y$ are i.i.d.
In view of
\begin{equation*}
{\bf 1}_{[y<x]} + {\bf 1}_{[y>x]} + {\bf 1}_{[y=x]} = {\bf 1}_{[ (x,y) \in {\mathbb R}^2]}
\end{equation*}
and (\ref{O12.2}), the left-hand sides of (\ref{O12.1.1}) and (\ref{O12.1.2}) are bounded by
\begin{equation*}
\tfrac 12 E_F \left( |\psi(X) \psi(Y)| \right) \leq \tfrac 12 E_F \left( \psi^2(X) \right).
\end{equation*}
The last inequality is valid by Cauchy-Schwarz, which also shows that equality can hold only if $\psi(X)$ is constant $F$-almost everywhere.
Indeed, for $\psi =1$ the left-hand sides of (\ref{O12.1.1}) and (\ref{O12.1.2}) equal \\
$P(Y<X) + \tfrac 12 P(Y=X) = P(X<Y) + \tfrac 12 P(Y=X)$ and the sum of these two terms equals \\
$P(Y<X) + P(Y>X) + P(Y=X) =1$.
\hfill
$\Box$

\begin{remark}\label{R0}
The continuous Opial inequality has also been generalized to fractional integral operators; for a recent comprehensive paper on this topic, see, e.g. \\ \cite{Vivas}.
For absolutely continuous distribution functions $F$ with density $f$ in (\ref{A}) and (\ref{B}) there is some similarity to fractional integral operators, in which $f(y)$ is replaced by a function of $(x-y)$, in the notation of (\ref{A}) and (\ref{B}).
\end{remark}

\begin{remark}\label{R1}
If $F$ has no point masses, i.e., if $P(X=Y) =0$ holds, then (\ref{O12.1.1}) simplifies to
\begin{equation*}
E_F \left( \left| E_F \left( \psi(Y) {\bf 1}_{[Y < X]} \mid X \right) \psi(X) \right| \right) \leq \tfrac 12 E_F\left( \psi^2(X) \right).
\end{equation*}
In particular, if $X$ and $Y$ have a uniform distribution on the interval $[0,h]$, then this reduces to Opial's inequality (\ref{A4})
\begin{equation}\label{O7}
\int_0^h \left| \left( \int_0^x \psi(y) \, dy \right) \psi(x) \right| \, dx \leq \frac h2 \int_0^h \psi^2(x) \, dx.
\end{equation}
In Opial's notation we have the absolutely continuous function $x(t) = \int_0^t \psi(y) \, dy$, which is Lebesgue almost everywhere differentiable with derivative $x'(t) = \psi(t)$ and $x(0)=0$.

With $x(t) =  \int_t^{2h} \psi(y) \, dy,\ h \leq t \leq 2h,$ we also have $x(2h) = 0$.
Then by symmetry (\ref{O7}) implies
\begin{equation*}
\int_h^{2h} \left| \left( \int_x^{2h} \psi(y) \, dy \right) \psi(x) \right| \, dx \leq \frac h2 \int_h^{2h} \psi^2(x) \, dx,
\end{equation*}
which combined with (\ref{O7}) itself yields inequality (\ref{A3}) of Opial, where $h$ is replaced by $2h$, namely
\begin{equation*}
\int_0^{2h} \left| \left( \int_0^x \psi(y) \, dy \right) \psi(x) \right| \, dx
\leq \frac h2 \int_0^h \psi^2(x) \, dx + \frac h2 \int_h^{2h} \psi^2(x) \, dx = \frac {2h}4 \int_0^{2h} \psi^2(x) \, dx.
\end{equation*}
Note that for continuity of $x$ at $h$ we need $ \int_0^h \psi(y) \, dy = \int_h^{2h} \psi(y) \, dy$.
In a similar way as above, two inequalities as in (\ref{O12.1.1}) may be glued together; see Corollary \ref{Cor} below.
\end{remark}

\begin{remark}\label{R2}
\cite{Lasota} presented a discrete Opial inequality, an analog of (\ref{A3}).
Its proof contains the more basic discrete version of (\ref{A4}) as formulated in Theorem 5.2.4 of \cite{Agarwal}.
With $a_i = \nabla x_i$ the backward shift operator, their Theorem 5.2.4 presents the inequality
\begin{equation}\label{O9.1}
\sum_{i=1}^N \left| \sum_{j=1}^i a_i a_j \right| \leq \frac{N+1}2 \sum_{i=1}^N a_i^2,\quad a_i \in {\mathbb R},\ i=1,\dots,N.
\end{equation}
With $F$ uniform on $\{ 1, \dots, N \}$ and $\ \psi(i)=a_i,\ i=1,\dots,N,$ the second inequality in (\ref{O12.1.1}) yields
\begin{equation}\label{O9.2}
\sum_{i=1}^N \sum_{j=1}^i \left| a_i a_j \right| \leq \frac{N+1}2 \sum_{i=1}^N a_i^2,\quad a_i \in {\mathbb R},\ i=1,\dots,N,
\end{equation}
which is equivalent to (\ref{O9.1}); replace $a_i$ by $|a_i|$.
\end{remark}

\begin{corollary}\label{Cor}
Let $X$ and $Y$ be i.i.d. with distribution function $F$ and let $c$ be a constant satisfying $P_F(X \leq c) = p,\, 0<p<1.$
Then the inequalities
\begin{eqnarray}\label{O10}
&& E_F \left( \left| E_F \left( \psi(Y) \left\{ {\bf 1}_{[Y<X]} + \tfrac 12 {\bf 1}_{[Y=X]} \right\} \mid X, Y \leq c \right) \psi(X) \right|
\mid X \leq c, Y \leq c \right) \nonumber \\
&& \hspace{2em} +\, E_F \left( \left| E_F \left( \psi(Y) \left\{ {\bf 1}_{[Y>X]} + \tfrac 12 {\bf 1}_{[Y=X]} \right\} \mid X, Y > c \right) \psi(X) \right| \mid X > c, Y > c \right) \nonumber \\
&& \leq E_F \left( \left| \psi(X) \psi(Y) \right| \left\{ {\bf 1}_{[Y<X]} + \tfrac 12 {\bf 1}_{[Y=X]} \right\} \mid X \leq c, Y \leq c \right) \\
&& \hspace{2em} +\, E_F \left( \left| \psi(X) \psi(Y) \right| \left\{ {\bf 1}_{[Y<X]} + \tfrac 12 {\bf 1}_{[Y=X]} \right\} \mid X > c, Y > c \right) \nonumber \\
&& \hspace{8em} \leq \tfrac 12 E_F\left( \psi^2(X) \left[ {\bf 1}_{[X \leq c]} / p + {\bf 1}_{[X > c]} / (1-p) \right] \right) \nonumber
\end{eqnarray}
hold with equalities if and only if $\psi$ is $F$-almost everywhere constant on $(-\infty,c]$ and on $(c,\infty)$ with possibly different constants. \\
In case $F$ is continuous and $c$ is a median of $F$, the inequalities become
\begin{eqnarray}\label{O10.1}
\lefteqn{ E_F \left( \left| E_F \left( \psi(Y) {\bf 1}_{[Y<X]} \mid X, Y \leq c \right) \psi(X) \right| {\bf 1}_{[X \leq c, Y \leq c]} \right) \nonumber } \\
&& \hspace{2em} +\, E_F \left( \left| E_F \left( \psi(Y) {\bf 1}_{[Y>X]} \mid X, Y > c \right) \psi(X) \right| {\bf 1}_{[X > c, Y > c]} \right) \\
&& \leq E_F \left( \left| \psi(X) \psi(Y) \right| {\bf 1}_{[Y<X]} {\bf 1}_{[X \leq c, Y \leq c]} \right)
   +\, E_F \left( \left| \psi(X) \psi(Y) \right| {\bf 1}_{[Y<X]} {\bf 1}_{[X > c, Y > c]} \right) \nonumber \\
&& \hspace{20em} \leq \tfrac 14 E_F\left( \psi^2(X) \right). \nonumber
\end{eqnarray}
\end{corollary}
\noindent
{\bf Proof} \\
Recall that $X$ and $Y$ are i.i.d.
Consequently they are i.i.d. under the condition $X>c, Y>c$ too.
Therefore the first inequality of (\ref{O10}) holds; cf. (\ref{O12.2}).
Applying Theorem \ref{Opial} to both expectations on the right-hand side of the first inequality
conditionally on $X \leq c, Y \leq c$ and $X > c, Y > c$, respectively,
we see that the left-hand sides of (\ref{O10}) are bounded from above by
\begin{eqnarray*}
\lefteqn{ \tfrac 12 E_F\left( \psi^2(X) \mid X \leq c,\, Y \leq c \right)
        + \tfrac 12 E_F\left( \psi^2(X) \mid X > c,\, Y > c \right) } \\
&& \hspace{5em} = \tfrac 12 E_F \left( \psi^2(X) \left[ {\bf 1}_{[X \leq c]} / p + {\bf 1}_{[X > c]} / (1-p) \right] \right). \nonumber
\end{eqnarray*}
As in Theorem \ref{Opial}, equalities hold if and only if $\psi$ is constant $F$-almost everywhere on $(-\infty,c]$ and on $(c,\infty)$ with possibly different constants.\\
For the case where $F$ is continuous and $c$ is a median of $F$, we have $P_F(Y=X)=0, \, p=1/2$,
\begin{eqnarray*}
\lefteqn{ E_F \left( \left| E_F \left( \psi(Y) {\bf 1}_{[Y<X]} \mid X, Y \leq c \right) \psi(X) \right| \mid X \leq c, Y \leq c \right) } \\
&& = 4 \times E_F \left( \left| E_F \left( \psi(Y) {\bf 1}_{[Y<X]} \mid X, Y \leq c \right) \psi(X) \right| {\bf 1}_{[X \leq c, Y \leq c]} \right)
\nonumber
\end{eqnarray*}
and similarly for the subsequent three terms.
\hfill
$\Box$

\begin{remark}\label{R3}
Taking $F$ uniform on $[0,h]$ in (\ref{O10.1}) we obtain $c=h/2,\, p=1/2$ and
\begin{eqnarray*}
\lefteqn{ \int_0^{h/2} \left| \psi(x) \int_0^x \psi(y) \, dy \right| \, dx + \int_{h/2}^h \left| \psi(x) \int_x^h \psi(y) \, dy \right| \, dx } \\
&& \leq  \int_0^{h/2} \int_0^x \left| \psi(x) \psi(y) \right| \, dy \, dx + \int_{h/2}^h \int_x^h \left| \psi(x) \psi(y) \right| \, dy \, dx
\leq \frac h4 \int_0^h \psi^2(x) \, dx.
\end{eqnarray*}
With $\int_x^h \psi = - \int_0^x \psi$ these inequalities yield
\begin{equation*}
\int_0^h \left| \int_0^x \psi(y) \, dy \, \psi(x) \right| \, dx \leq \frac h4 \int_0^h \psi^2(x) \, dx,
\end{equation*}
our equivalent of (\ref{A4}); see also Remark \ref{R1}.
\end{remark}

\begin{remark}\label{R4}
With $F$ uniform on $\{1,\dots,K,\dots,N\},\ c=K,\ p=K/N$ and $\psi(i)=a_i,\, i=1, \dots,N,$ inequality (\ref{O10}) becomes
\begin{eqnarray*}
\lefteqn{ \frac 1{K^2} \sum_{i=1}^K \left| a_i \left\{ \sum_{j=1}^i a_j - \tfrac 12 a_i \right\} \right|
+ \frac 1{(N-K)^2} \sum_{i=K+1}^N \left| a_i \left\{ \sum_{j=i+1}^N a_j + \tfrac 12 a_i \right\} \right| \nonumber } \\
&& \hspace{15em} \leq \frac 1{2K} \sum_{i=1}^K a_i^2 + \frac 1{2(N-K)} \sum_{i=K+1}^N a_i^2,
\end{eqnarray*}
which for $N=2K$ and $\sum_{i=1}^N \psi(i) = \sum_{i=1}^N a_i = 0$ reduces to
\begin{equation}\label{O15}
\sum_{i=1}^N \left| a_i \left\{ \sum_{j=i+1}^N a_j + \tfrac 12 a_i \right\} \right|
= \sum_{i=1}^N \left| a_i \left\{ \sum_{j=1}^{i-1} a_j + \tfrac 12 a_i \right\} \right|
\leq \frac N4 \sum_{i=1}^N a_i^2, \quad \left( \sum_{i=1}^N a_i = 0 \right).
\end{equation}
This inequality is sharp as can be seen via the choices $a_1 = \dots = a_K =1,\, a_{K+1} = \dots = a_N = -1$.
For odd $N$ this inequality also holds as we will show next.
Define $b_h = \psi(\lceil h/2 \rceil) = a_{\lceil h/2 \rceil}, \, h = 1, \dots, 2N.$
In view of $\sum_{h=1}^{2N} b_h = 2 \sum_{i=1}^N a_i =0$ inequality (\ref{O15}) yields
\begin{equation*}
\sum_{h=1}^{2N} \left| b_h \left\{ \sum_{k=1}^{h-1} b_k + \tfrac 12 b_h \right\} \right|
\leq \frac {2N}4 \sum_{h=1}^{2N} b_h^2 = N \sum_{i=1}^N a_i^2.
\end{equation*}
Discerning between odd and even $h$, we see that this inequality can be rewritten in terms of the $a_i$ as
\begin{equation*}
\sum_{i=1}^N \left( \frac 12 \left| a_i \left\{ \sum_{j=1}^{i-1} a_j + \tfrac 34 a_i \right\} \right|
+ \, \frac 12 \left| a_i \left\{ \sum_{j=1}^{i-1} a_j + \tfrac 14 a_i \right\} \right| \right) \leq \frac N4 \sum_{i=1}^N a_i^2.
\end{equation*}
Since the absolute value function $| \cdot|$ is convex this inequality implies (\ref{O15}) for odd $N$.

The discrete Opial inequality in Theorem 1.2 of \cite{Lasota} states
\begin{equation}\label{O18}
\sum_{i=1}^N \left| a_i \sum_{j=1}^{i-1} a_j \right|
\leq \frac 12 \left\lfloor \frac {N+1}2 \right\rfloor \sum_{i=1}^N a_i^2,
\quad \left( \sum_{i=1}^N a_i = 0 \right);
\end{equation}
see also Theorem 5.2.1 of \cite{Agarwal}.
Again, for even $N=2K$ equality holds for $a_1 = \dots = a_K =1,\, a_{K+1} = \dots = a_N = -1$. \\
Note that the left-hand side of Lasota's inequality (\ref{O18}),
\begin{equation*}
\sum_{i=1}^N \left| a_i \sum_{j=1}^{i-1} a_j \right| = \sum_{i=1}^N \left| a_i \sum_{j=i}^N a_j \right|,
\end{equation*}
lacks the symmetry that the left-hand side of our inequality (\ref{O15}) has.
\end{remark}

\section{A generalization of Opial's inequality for the $n$-th derivative}

Let $X$ be a random variable with distribution function $F$, let $\psi$ be a measurable function and let $n$ be a natural number.
We define the $n$-th order integral
\begin{equation*}
I_{F,n,\psi}(x) = \int_{(-\infty,x)} \int_{(-\infty, x_{n-1})} \cdots \int_{(-\infty, x_1)} \psi(x_0)\, dF(x_0) \cdots dF(x_{n-2}) \, dF(x_{n-1}).
\end{equation*}
Note that for $n=1$, we have
\begin{equation*}
I_{F,1,\psi}(x) = \int_{(-\infty,x)} \psi(x_0) \, dF(x_0) = E_F \left( \psi(Y) {\bf 1}_{[Y < X]} \mid X = x \right).
\end{equation*}
We generalize the $n$-th order derivative version of Opial's inequality, as given by \cite{Das}, and improve the constant in this generalization to its optimal value.
\begin{theorem}\label{Das}
With the above notation
\begin{equation}\label{M3}
E_F \left( \left| I_{F,n,\psi}(X) \, \psi(X) \right| \right) \leq \frac 1 {(n+1)!} \ E_F \left( \psi^2(X) \right)
\end{equation}
holds with equality if and only if $F$ is continuous and $\psi$ is a constant $F$-almost everywhere.
\end{theorem}
\noindent
{\bf Proof} \\
Without loss of generality we assume $E_F\left( \psi^2(X) \right) < \infty$.
In terms of integrals, the left-hand side of (\ref{M3}) may be written as
\begin{eqnarray*}
\lefteqn{\int_{-\infty}^\infty \left| \int_{(-\infty,x_n)} \int_{(-\infty, x_{n-1})} \cdots \int_{(-\infty, x_1)}
\psi(x_0)\, \psi(x_n) \, dF(x_0) \cdots dF(x_{n-2}) \, dF(x_{n-1}) \right| } \\
&& \hspace{32em} dF(x_n),
\end{eqnarray*}
which is bounded from above by
\begin{eqnarray}\label{M4.1}
\lefteqn{ \int_{-\infty}^\infty \int_{(-\infty,x_n)} \int_{(-\infty, x_{n-1})} \cdots \int_{(-\infty, x_1)}
\left| \psi(x_0)\, \psi(x_n) \right| } \\
&& \hspace{18em} dF(x_0) \cdots dF(x_{n-2}) \, dF(x_{n-1}) \, dF(x_n). \nonumber
\end{eqnarray}
Let $S_{n+1}$ be the set of permutations of $\{0,1,\dots,n\}$.
Consider the disjoint union
\begin{equation}\label{M5}
U_{n+1} = \bigcup_{\pi \in S_{n+1}} \left\{ (x_0,\dots, x_n) \in {\mathbb R}^{n+1} \mid x_{\pi(0)} < x_{\pi(1)} , \cdots < x_{\pi(n)} \right\}
\subsetneqq {\mathbb R}^{n+1}.
\end{equation}
Since each of the $(n+1)!$ permutations $\pi \in S_{n+1}$ applied to (\ref{M4.1}) yields the same value, we obtain
\begin{eqnarray}\label{M6}
\lefteqn{ (n+1)! \ E_F \left( \left| I_{F,n,\psi}(X) \, \psi(X) \right| \right) } \\
&& \leq \int_{-\infty}^\infty \int_{-\infty}^\infty \int_{-\infty}^\infty \cdots \int_{-\infty}^\infty
\left| \psi(x_0)\, \psi(x_n) \right| \, dF(x_0) \cdots dF(x_{n-2}) \, dF(x_{n-1}) \, dF(x_n) \nonumber \\
&& = \left( E_F \left( \left| \psi(X) \right| \right) \right)^2 \leq  E_F \left( \psi^2(X) \right). \nonumber
\end{eqnarray}
In view of the inequality sign in (\ref{M5}) the first inequality in (\ref{M6}) can be an equality if and only if $F$ has no point masses, i.e., if and only if $F$ is continuous.
Because of the Cauchy-Schwarz inequality, the second inequality is an equality if and only if $\psi$ is constant $F$-almost everywhere.
\hfill
$\Box$

\begin{remark}\label{Das2}
Let $F$ be the uniform distribution function on $[a,b]$ as in Theorem 1 of \cite{Das}.
Then we have
\begin{equation*}
I_{F,n,\psi}(x) = (b-a)^{-n} \int_a^x \int_a^{x_{n-1}} \cdots \int_a^{x_1} \psi(x_0)\, dx_0 \cdots dx_{n-2} \, dx_{n-1}
\end{equation*}
and we see that the $n$-th derivative of $(b-a)^n \, I_{F,n,\psi}(x)$ equals $\psi(x)$.
We conclude that our Theorem \ref{Das} generalizes Theorem 1 of \cite{Das} with its constant $K = \sqrt{n/(2n-1)}\,/(2(n!))$ replaced by the optimal $1/(n+1)!$. (We have $K \geq 1/(n+1)!$ with equality only if $n=1$.)
\end{remark}

\begin{remark}\label{Das3}
For $n=1$ and $F$ continuous Theorem \ref{Das} yields (\ref{O12.1.1}) of Theorem \ref{Opial}.
\end{remark}

\section{A sharp inequality for $n=2$ and distributions with atoms}\label{n2}

If $F$ has point masses, then the probability $P_F(X_0=X_1)$ is positive and hence $P\left( X \in {\mathbb R}^{n+1} \setminus U_{n+1} \right)$ is.
Therefore, in order to improve inequality (\ref{M3}) to a sharp inequality for distributions with point masses, we have to study the gap between $U_{n+1}$ and ${\mathbb R}^{n+1}$; recall the inequality sign in (\ref{M5}).
We do this for $n=2$ as follows.
Let
\begin{eqnarray*}
&& U_3 = \bigcup_{\pi \in S_3} \left\{ (x_0, x_1, x_2) \in {\mathbb R}^3 \mid x_{\pi(0)} < x_{\pi(1)} < x_{\pi(2)} \right\}, \quad {\rm disjoint~union~of~6~sets}, \nonumber \\
&& V_1 = \bigcup_{\pi \in S_3} \left\{ (x_0, x_1, x_2) \in {\mathbb R}^3 \mid x_{\pi(0)} = x_{\pi(1)} < x_{\pi(2)} \right\}, \quad {\rm disjoint~union~of~3~sets}, \\
&& V_2 = \bigcup_{\pi \in S_3} \left\{ (x_0, x_1, x_2) \in {\mathbb R}^3 \mid x_{\pi(0)} < x_{\pi(1)} = x_{\pi(2)} \right\}, \quad {\rm disjoint~union~of~3~sets}, \nonumber \\
&& W = \left\{ (x_0, x_1, x_2) \in {\mathbb R}^3 \mid x_0 = x_1 = x_2 \right\}. \nonumber
\end{eqnarray*}
Now
\begin{equation*}
{\mathbb R}^3 = U_3 \cup V_1 \cup V_2 \cup W
\end{equation*}
is a partition of ${\mathbb R}^3$ consisting of 13 sets.
By analogy with (\ref{M6}), using the Tonelli-Fubini theorem, we have
\begin{eqnarray}\label{two3}
\lefteqn{  E_F \left( \psi^2(X) \right) \geq \left( E_F \left( \left| \psi(X) \right| \right) \right)^2 \nonumber } \\
&& = \int_{-\infty}^\infty \int_{-\infty}^\infty \int_{-\infty}^\infty
\left| \psi(x_0)\, \psi(x_2) \right| \left[ {\bf 1}_{U_3} + {\bf 1}_{V_1} + {\bf 1}_{V_2} + {\bf 1}_W \right] \, dF(x_0) \, dF(x_1) \, dF(x_2) \nonumber \\
&& = 6 \times \int_{-\infty}^\infty \int_{(-\infty, x_2)} \int_{(-\infty, x_1)}
\left| \psi(x_0)\, \psi(x_2) \right| \, dF(x_0) \, dF(x_1) \, dF(x_2) \nonumber \\
&& \hspace{2em} + \, 3 \times \iiint_{x_0=x_1<x_2} \left| \psi(x_0)\, \psi(x_2) \right| \, dF(x_0) \, dF(x_1) \, dF(x_2) \nonumber \\
&& \hspace{2em} + \, 3 \times \iiint_{x_0<x_1=x_2} \left| \psi(x_0)\, \psi(x_2) \right| \, dF(x_0) \, dF(x_1) \, dF(x_2) \\
&& \hspace{2em} + \, \iiint_{x_0=x_1=x_2} \left| \psi(x_0)\, \psi(x_2) \right| \, dF(x_0) \, dF(x_1) \, dF(x_2) \nonumber \\
&& = 6 \times \int_{-\infty}^\infty \int_{(-\infty, x_2)} \int_{(-\infty, x_1)}
\left| \psi(x_0) \right| \, dF(x_0) \, dF(x_1) \, \left| \psi(x_2) \right| dF(x_2) \nonumber \\
&& \hspace{2em} + \, 3 \times \int_{-\infty}^\infty \int_{(-\infty, x_2)} \left| \psi(x_0) \right| P_F(X=x_0)\, dF(x_0) \left|\psi(x_2) \right| \, dF(x_2) \nonumber \\
&& \hspace{2em} + \, 3 \times \int_{-\infty}^\infty \int_{(-\infty, x_2)} \left| \psi(x_0) \right| \, dF(x_0) \, P_F(X=x_2) \left|\psi(x_2) \right| \, dF(x_2) \nonumber \\
&& \hspace{2em} + \, \int_{-\infty}^\infty \psi^2(x_2)\, P_F^2(X=x_2) \, dF(x_2). \nonumber
\end{eqnarray}

With the notation
\begin{eqnarray*}
&& J_{F,\psi}(x) = \int_{(-\infty,x)} \int_{(-\infty, y)} \left| \psi(z) \right| \, dF(z) \, dF(y), \nonumber \\
&& J_{F,\psi, {\rm D}}(x) = \int_{(-\infty, x)} \left| \psi(z) \right| \left[ P_F(X=z) + P_F(X=x) \right] \, dF(z), \\
&& p_F(x) = P_F (X=x), \quad x \in {\mathbb R}, \nonumber
\end{eqnarray*}
we have proved the inequality
\begin{theorem}\label{n=2}
Let $X$ be a random variable with distribution function $F$ on $\mathbb R$, and let $\psi : {\mathbb R} \to {\mathbb R}$ be a measurable function. The inequality
\begin{equation}\label{two4}
6 \times E_F \left( J_{F, \psi}(X) \left| \psi(X) \right| \right) + 3 \times E_F \left( J_{F, \psi, {\rm D}}(X) \left| \psi(X) \right| \right)
\leq E_F \left( \psi^2(X) \left[ 1 - p_F^2(X) \right] \right)
\end{equation}
holds with equality if and only if $\psi$ is constant $F$-almost everywhere.
\end{theorem}
\noindent
{\bf Proof} \\
The only inequality in (\ref{two3}) stems from Cauchy-Schwarz, which means that equality in (\ref{two4}) holds if and only $\psi$ is constant $F$-almost everywhere.
\hfill
$\Box$ \\

If $F$ is continuous, i.e. if $F$ has no point masses, then inequality (\ref{two4}) reduces to inequality (\ref{M3}) with $n=2$.
\begin{remark}\label{Rtwo}
Let $F$ be uniform on $\{1,\dots,N\}$ and denote $|\psi(i)| = a_i,\ i=1, \dots, N$.
Some computation shows that inequality (\ref{two4}) may be rewritten then as
\begin{equation*}
6 \, \sum_{i=1}^N \sum_{j=1}^i (i-j) a_i a_j \leq \left( N^2 -1 \right) \sum_{i=1}^N a_i^2,
\end{equation*}
which for $a_i,\, i=1,\dots,N,$ constant is an equality indeed.
\end{remark}

\section{Opial inequalities with weights}

Applying our approach of Theorem \ref{Opial} we derive Opial inequalities with a weight function.
\begin{theorem}\label{WeightedOpial}
Let $X$ and $Y$ be independent and identically distributed (i.i.d.) random variables with distribution function $F$ on $\mathbb R$, let $\psi : {\mathbb R} \to {\mathbb R}$ be a measurable function, and let $\chi : {\mathbb R} \to {\mathbb R}$ be a nonnegative weight function.
The inequalities
\begin{eqnarray}\label{weight.1}
\lefteqn{ E_F \left( \left| E_F \left( \psi(Y) \left\{ {\bf 1}_{[Y<X]} + \tfrac 12 {\bf 1}_{[Y=X]} \right\} \mid X \right) \psi(X) \right| \chi(X) \right) \nonumber } \\
&& \leq E_F \left( |\psi(X)\psi(Y)| \chi(X) \left\{ {\bf 1}_{[Y<X]} + \tfrac 12 {\bf 1}_{[Y=X]} \right\} \right) \\
&& \leq \tfrac 12 E_F \left( \psi^2(X)
\left[ \chi(X) \left\{ {\bf 1}_{[Y < X]} + \tfrac 12 {\bf 1}_{[Y=X]} \right\} + \chi(Y) \left\{ {\bf 1}_{[Y > X]} + \tfrac 12 {\bf 1}_{[Y=X]} \right\} \right] \right) \nonumber
\end{eqnarray}
and
\begin{eqnarray}\label{weight.2}
\lefteqn{ E_F \left( \left| E_F \left( \psi(Y) \left\{ {\bf 1}_{[Y>X]} + \tfrac 12 {\bf 1}_{[Y=X]} \right\} \mid X \right) \psi(X) \right| \chi(X) \right) \nonumber } \\
&& \leq E_F \left( |\psi(X)\psi(Y)| \chi(X) \left\{ {\bf 1}_{[Y>X]} + \tfrac 12 {\bf 1}_{[Y=X]} \right\} \right) \\
&& \leq \tfrac 12 E_F \left( \psi^2(X)
\left[ \chi(X) \left\{ {\bf 1}_{[Y > X]} + \tfrac 12 {\bf 1}_{[Y=X]} \right\} + \chi(Y) \left\{ {\bf 1}_{[Y < X]} + \tfrac 12 {\bf 1}_{[Y=X]} \right\} \right] \right) \nonumber
\end{eqnarray}
hold with equality if $\psi$ is constant $F$-almost everywhere.
In case $\chi$ is nonincreasing or nondecreasing, the upperbounds in (\ref{weight.1}) and (\ref{weight.2}), respectively, are bounded from above by $\tfrac 12  E_F \left( \psi^2(X) \chi(X) \right)$.
\end{theorem}
\noindent
{\bf Proof} \\
We focus on the second inequality in (\ref{weight.1}).
Since $X$ and $Y$ are i.i.d. we have
\begin{eqnarray*}
\lefteqn{ E_F \left( |\psi(X)\psi(Y)| \chi(X) \left\{ {\bf 1}_{[Y<X]} + \tfrac 12 {\bf 1}_{[Y=X]} \right\} \right) } \\
&& = E_F \left( |\psi(X)\psi(Y)| \chi(Y) \left\{ {\bf 1}_{[Y>X]} + \tfrac 12 {\bf 1}_{[Y=X]} \right\} \right) \\
&& = \frac 12 E_F \left( |\psi(X)\psi(Y)| \left[ \chi(X) \left\{ {\bf 1}_{[Y<X]} + \tfrac 12 {\bf 1}_{[Y=X]} \right\}
                        + \chi(Y) \left\{ {\bf 1}_{[Y>X]} + \tfrac 12 {\bf 1}_{[Y=X]} \right\} \right] \right) \\
&& \leq \frac 12 E_F \left( \psi^2(X) \left[ \chi(X) \left\{ {\bf 1}_{[Y<X]} + \tfrac 12 {\bf 1}_{[Y=X]} \right\}
                        + \chi(Y) \left\{ {\bf 1}_{[Y>X]} + \tfrac 12 {\bf 1}_{[Y=X]} \right\} \right] \right).
\end{eqnarray*}
Again, since $X$ and $Y$ are i.i.d. the left- and right-hand side of (\ref{weight.1}) are equal for $\psi=1$. \\
In case $\chi$ is nonincreasing $\chi(Y) {\bf 1}_{[Y>X]} \leq \chi(X) {\bf 1}_{[Y>X]}$ holds, which implies the last statement of the theorem for (\ref{weight.1}).

The proof for (\ref{weight.2}) is analogous.
\hfill
$\Box$ \\

\begin{remark}\label{W1}
Note that inequalities (\ref{weight.1}) and (\ref{weight.2}) generalize inequalities (\ref{O12.1.1}) and (\ref{O12.1.2}), as can be seen by choosing $\chi=1$.
\end{remark}

\begin{remark}
With $F$ uniform on the interval $(0,h)$, the inequalities from Theorem \ref{WeightedOpial} yield
\begin{eqnarray}\label{weight.4}
&& \int_0^h \left| \int_0^x \psi(y)\, dy \, \psi(x) \right| \chi(x) \, dx
    \leq \frac 12 \int_0^h \psi^2(x) \left\{ \chi(x) x + \int_x^h \chi(y) \, dy \right\} \, dx, \\
&& \int_0^h \left| \int_x^h \psi(y)\, dy \, \psi(x) \right| \chi(x) \, dx
    \leq \frac 12 \int_0^h \psi^2(x) \left\{ \chi(x) (h-x) + \int_0^x \chi(y) \, dy \right\} \, dx. \nonumber
\end{eqnarray}
\cite{Troy} considers the special case of $\chi(x) = x^p,\ p>-1,$ and states: "It remains an open problem to determine the sharpness of ... and (5).", which in our notation is
\begin{equation}\label{weight.5}
\int_0^h x^p \left| \int_0^x \psi(y)\, dy \, \psi(x) \right| \, dx \leq \frac {h^{p+1}}{2\sqrt{p+1}} \int_0^h \psi^2(x) \, dx.
\end{equation}
However, equalities in (7) and (10) of \cite{Troy} imply that $\psi(x)$ has to be a multiple of $x$ in order for this inequality to be an equality.
With $\psi(x) = x$, inequality (\ref{weight.5}) becomes $h^{p+4}/(2(p+4)) \leq h^{p+4}/(6\sqrt{p+1})$, which is a strict inequality for all $p$.
With $\chi(x) = x^p,\ p>-1,$ the first one of our inequalities in (\ref{weight.4}) yields
\begin{equation*}
\int_0^h x^p \left| \int_0^x \psi(y)\, dy \, \psi(x) \right| \, dx
\leq \frac {h^{p+1}}{2(p+1)} \int_0^h \psi^2(x) \left\{1 + p\left( \frac xh \right)^{p+1} \right\} \, dx,
\end{equation*}
which is less elegant than (\ref{weight.5}), but it is sharp.
\end{remark}

\begin{remark}
With $F$ uniform on $\{ 1, \dots, N \}$ and $\ \psi(i)=a_i, \chi(i)=w_i,\ i=1,\dots,N,$ the second inequality in (\ref{weight.1}) yields
\begin{equation*}
\sum_{i=1}^N \sum_{j=1}^i \left| a_i a_j \right| w_i \leq \frac 12 \sum_{i=1}^N a_i^2 \left[(i+1)w_i + \sum_{j=i+1}^N w_j \right],
\end{equation*}
which generalizes (\ref{O9.2}).
\end{remark}

\section{The Wirtinger inequality}

Closely related to Opial's inequality is Wirtinger's inequality (\ref{A2}), which is an equality
if and only if $x(t)$ is a constant times $\sin(\pi x/h)$.
We apply our approach with integration with respect to a distribution function and get the following result.
\begin{theorem}\label{Wirtinger}
{\bf Wirtinger's inequality} \\
Let $X$ and $Y$ be independent and identically distributed random variables with absolutely continuous distribution function $F$ on $\mathbb R$, and let $\psi : {\mathbb R} \to {\mathbb R}$ be a measurable function.
With $E_F(\psi(X))=0$ the inequality
\begin{equation}\label{Wirt2}
E_F \left( \left[ E_F \left(\psi(Y) {\bf 1}_{[Y<X]} \mid X \right) \right]^2 \right) \leq \frac 1{\pi^2} \, E_F \left( \psi^2(X) \right)
\end{equation}
holds with equality if and only if $\psi$ is $F$-almost everywhere a constant times $\cos ( \pi F )$ .
\end{theorem}
\noindent
{\bf Proof} \\
We assume $E_F(\psi^2(X)) = \int_{-\infty}^\infty \psi^2(x) \, dF(x) = \int_{-\infty}^\infty \psi^2(x) f(x) \, dx < \infty$,
where $f$ is a density of $F$.
In view of
\begin{eqnarray*}
\lefteqn{ \frac d{dx} \left( \left( \int_{-\infty}^x \psi(y) f(y) \, dy \right)^2 \cot (\pi F(x)) \right) \nonumber }\\
&& = 2 \int_{-\infty}^x \psi(y) f(y) \, dy \ \psi(x) f(x) \cot(\pi F(x)) \\
&& \hspace{1em} - \left( \int_{-\infty}^x \psi(y) f(y) \, dy \right)^2
\left[ 1 + \cot^2(\pi F(x)) \right] \pi f(x), \quad F-a.e., \nonumber
\end{eqnarray*}
a similar equation holds as in a classic proof of Wirtinger's inequality, namely
\begin{eqnarray}\label{Wirt4}
\lefteqn{ \left( \int_{-\infty}^x \psi(y) f(y) \, dy \right)^2 f(x)
        + \left( \frac {\psi(x)}{\pi} - \left( \int_{-\infty}^x \psi(y) f(y) \, dy \right) \cot (\pi F(x)) \right)^2 f(x) \nonumber } \\
&&  = \frac{\psi^2(x)}{\pi^2} f(x) - \frac 2\pi \psi(x) \int_{-\infty}^x \psi(y) f(y) \, dy \, \cot (\pi F(x)) \ f(x) \\
&& \hspace{1em} + \left( \int_{-\infty}^x \psi(y) f(y) \, dy \right)^2 \cot^2 (\pi F(x)) \ f(x)
        + \left( \int_{-\infty}^x \psi(y) f(y) \, dy \right)^2 f(x) \nonumber \\
&& = \frac{\psi^2(x)}{\pi^2} f(x) - \frac 1\pi \, \frac d{dx} \left( \left( \int_{-\infty}^x \psi(y) f(y) \, dy \right)^2 \cot (\pi F(x)) \right). \nonumber
\end{eqnarray}
Because of $\cot y < 1/y$ for $0<y$ small, we have
\begin{eqnarray*}
\lefteqn{ \limsup_{x \to -\infty} \left( \int_{-\infty}^x \psi(y) f(y) \, dy \right)^2 \cot (\pi F(x)) } \\
&& \leq \limsup_{x \to -\infty} \left( \int_{-\infty}^x \psi^2(y) f(y) \, dy \right) F(x) \, \frac 1{\pi F(x)} =0
\end{eqnarray*}
and analogously for $x \to \infty$.
This implies that the integral with respect to $x$ of the last term in (\ref{Wirt4}) vanishes.
Furthermore the integral of the second term in (\ref{Wirt4}) is clearly nonnegative and equals 0 if and only if $\psi$ is $F$-almost everywhere a constant times $\cos ( \pi F )$, as some computation shows.
Integration of (\ref{Wirt4}) proves the theorem.
\hfill
$\Box$ \\

\begin{remark}
With $F$ uniform on $[0,h]$, i.e., with $f(x) = 1/h \, {\bf 1}_{[0 \leq x \leq h]}$, and with $x(t) = \int_0^t \psi(y) \, dy$, inequality (\ref{Wirt2}) yields (\ref{A2}).
\end{remark}

\begin{remark}
Theorem \ref{Wirtinger} contains a known Wirtinger inequality with a weight.
In the notation of integrals inequality (\ref{Wirt2}) reads as
\begin{equation*}
\int_{-\infty}^\infty \left[ \int_{-\infty}^x \psi(y)f(y) \, dy \right]^2 f(x) \, dx
\leq \frac 1{\pi^2} \int_{-\infty}^\infty \psi^2(x) f(x) \, dx, \quad \int_{-\infty}^\infty \psi(x) f(x) \, dx = 0.
\end{equation*}
Writing $\chi = \psi f$, this inequality becomes
\begin{equation*}
\int_{-\infty}^\infty \left[ \int_{-\infty}^x \chi(y) \, dy \right]^2 f(x) \, dx
\leq \frac 1{\pi^2} \int_{-\infty}^\infty \chi^2(x) \frac 1{f(x)} \, dx, \quad \int_{-\infty}^\infty \chi(x) \, dx = 0,
\end{equation*}
or in the notation of Section 1 with $x(t) = \int _{-\infty}^t \chi(y) \, dy, \ w(t) = Wf(t),$
\begin{equation*}
\int_{-\infty}^\infty x^2(t)w(t) \, dt \leq \frac {W^2}{\pi^2} \int_{-\infty}^\infty (x'(t))^2 \frac 1{w(t)} \, dt,\quad x(-\infty)=x(\infty) = 0,
\end{equation*}
where $\chi$ and $x'$, respectively, vanish if $f$ does; cf. (\ref{A2}). \\
This generalizes Theorem 1.1 with $p=q=2$ of \cite{Giova} and (2.1) of Theorem 1 with $p=2$ of \cite{Dragomir}.
\end{remark}

\section{Conclusion}

There exist several versions of the Opial inequality.
We discern two classes, namely the class of two-sided inequalities, as we will call them, with $\int \psi \, dF = 0$ and the one-sided ones without this condition.
In Theorem \ref{Opial} we have presented a generalization of the classic one-sided inequality that contains both the continuous and the discrete inequalities with the optimal constants.
Our two-sided general inequality from Corollary \ref{Cor} contains the classic continuous Opial inequality with the optimal constant, but the discrete version of it deviates a little from the one by \cite{Lasota}.
However, our discrete inequality shows symmetry and Lasota's does not.

Theorem \ref{Das} presents a general Opial inequality for derivatives of order $n$.
For the continuous case it has the optimal constant, not known in the literature yet, as far as we can tell.
The case with $F$ having discontinuities, in particular the discrete case, is much more tricky, since one then has to partition ${\mathbb R}^{n+1}$ in a complicated way.
Moreover, the resulting inequality will then contain derivatives/differences of several orders.
In Section \ref{n2} this is illustrated for the case $n=2$.

Furthermore, Theorem \ref{WeightedOpial} generalizes our Opial inequalities from Theorem \ref{Opial} by introducing a weight function.
The inequalities obtained are sharp, but less elegant for a uniform distribution on an interval than the inequalities of \cite{Troy}, which are not sharp though.

Finally we presented our generalization of Wirtinger's inequality.
It automatically includes a Wirtinger inequality with a weight function.

These results lead us to the conclusion that it might be worthwhile to study inequalities involving (Lebesgue) integrals of functions and their derivatives by our technique, i.e., by integration with respect to distribution functions.
This approach is unifying in that the resulting inequalities often contain a discrete version of the inequality concerned as well.
Typically, the continuous and discrete versions are obtained by choosing uniform distribution functions on either an interval (of reals) or on an integer interval.

This unifying, generalizing approach was applied before in \cite{Wellner} and \\
\cite{Klaassen}.

{}
\end{document}